\documentclass{birkjour}
\usepackage{amsthm}
\usepackage{mathtools}
\usepackage{url}
\usepackage{helvet}
\usepackage{times}
\usepackage{amsmath}
\usepackage{amsfonts}
\usepackage{amssymb}
\usepackage{mathrsfs}
\usepackage{bbm}
\usepackage{floatflt}
\usepackage{graphics}
 \newtheorem{thm}{Theorem}[section]

 \theoremstyle{definition}
 
 \theoremstyle{remark}

 \numberwithin{equation}{section}
\def\triple{{\frak g}_{\frak t}}
\begin{document}
\title{Weil algebra, 3-Lie algebra and B.R.S. algebra}
\author{Viktor Abramov}
\address{Institute of Mathematics, University of Tartu\\
Liivi 2 -- 602, Tartu 50409, Estonia}
\email{viktor.abramov@ut.ee}


\date{}

\begin{abstract}
We consider the 3-Lie algebra induced by a Lie algebra with the help of an analog of a trace. We propose the extension of the Weil algebra of a Lie algebra to the Weil algebra of induced 3-Lie algebra by introducing in addition to an analog of connection and its curvature new elements and defining their differential by means of structure constants of 3-Lie algebra. We also propose a new approach to the universal B.R.S. algebra based on the quantum triple Nambu bracket.
\end{abstract}

\maketitle
\section{Introduction}
$n$-Lie algebra, where integer $n$ satisfies $n\geq 2$, can be regarded as a generalization of the concept of binary Lie algebra, which is the particular case of $n$-Lie algebra when $n=2$. An integer $n$ shows the number of elements of algebra necessary to compose a Lie bracket. Thus, a Lie bracket of $n$-Lie algebra $\frak g$ is an $n$-ary multilinear mapping $\frak g\times\frak g\times\ldots\times\frak g\;(n\; \mbox{times})\to\frak g$, which is totally skew-symmetric and satisfies the Filippov-Jacobi identity (also called fundamental identity).

\vskip.3cm
\noindent
A notion of $n$-Lie algebra was introduced and studied by V.T. Filippov in \cite{Filippov} and Y. Nambu in \cite{Nambu}, who proposed a generalization of Hamiltonian mechanics based on a triple (or, more generally, $n$-ary) bracket of functions determined on a phase space. This $n$-ary bracket can be regarded as an analog of the Poisson bracket in Hamiltonian mechanics and later this $n$-ary bracket was called the Nambu bracket. The dynamics of this generalization of Hamiltonian mechanics is described by the Nambu-Hamilton equation of motion that contains $n-1$ Hamiltonians. It should be pointed out that Y. Nambu proposed and developed his generalization of Hamiltonian mechanics with the goal to apply this approach to quarks model, where baryons are combinations of three quarks. Later it was shown that the Nambu bracket satisfied the Filippov-Jacobi identity and thus, the space of functions endowed with the $n$-ary Nambu bracket can be regarded as the example of $n$-Lie algebra. This direction of research has received further development in the paper \cite{Takhtajan}, where the author introduced a notion of Nambu-Poisson manifold, which can be regarded as an analog of the notion of Poisson manifold in Hamiltonian mechanics. 

\vskip.3cm
\noindent
The problem of quantization of Nambu bracket was considered in a number of papers and the first step in this direction was taken already by Y. Nambu, but so far this is the outstanding problem. In the paper \cite{Awata-Li-Minic-Yaneya} the authors propose the realization of quantum Nambu bracket defined by
\begin{equation}
[X,Y,Z]=\mbox{Tr}\,X\;[Y,Z]+\mbox{Tr}\,Y\;[Z,X]+\mbox{Tr}\,Z\;[X,Y],
\label{formula 1 in introduction}
\end{equation}
where $X,Y,Z$ are either square matrices or cubic matrices, and prove that this quantum Nambu bracket satisfies the Filippov-Jacobi identity. Hence a matrix Lie algebra endowed with the triple quantum Nambu bracket (\ref{formula 1 in introduction}) becomes the 3-Lie algebra which will be referred to as the induced 3-Lie algebra. The triple quantum Nambu bracket can be generalized by means of vector fields and differential forms on a smooth manifold \cite{AbramovIJGMMP} and by means of a cochain of the Eilenberg-Cartan complex of a Lie algebra \cite{Abramov2017Gent}.

\vskip.3cm
\noindent
Given a Lie algebra $\frak g$ one can associate to it the Weil algebra $W(\frak g)$. The Weil algebra of a Lie algebra is the tensor product of exterior algebra of dual space $\frak g^\ast$ and the symmetric algebra of the dual space $\frak g^\ast$. The Weil algebra is closely related to equivariant differential forms on a principal fiber bundle and it can be regarded as an algebraic structure for connections and their curvatures on a principal fiber bundle \cite{Quillen-Mathai}. It is also worth to mention that the Weil algebra is closely related to BRST algebra \cite{Dubois-Violette_Talon_Vialett}.

\vskip.3cm
\noindent
The aim of this paper is to extend the Weil algebra of a Lie algebra $\frak g$ to an analogous algebra of the induced 3-Lie algebra. The quantum Nambu bracket of induced 3-Lie algebra is defined by
$$
[x,y,z]=\omega(x)\,[y,z]+\omega(y)\,[z,x]+\omega(z)\,[x,y],
$$
where $\omega\in\frak g^\ast$ and $\omega([x,y])=0$.
In order to do this we introduce in addition to two elements of degree 1 and 2 of the Weil algebra new elements constructed by means of element $\omega$ of the dual space $\frak g^\ast.$ We extend the differential of the Weil algebra to these new elements and show that the property $d^2=0$ of the differential holds for these new elements. 

\section{3-Lie algebra induced by a Lie algebra}
An $n$-Lie algebra, where integer $n$ is greater or equal to 2, is a vector space equipped with an $n$-ary Lie bracket
$$
(x_1,x_2,\ldots,x_n)\in \frak g\times\frak g\times\ldots\times\frak g\;(n\;\mbox{times})\to [x_1,x_2,\ldots,x_n]\in \frak g,
$$
which is totally skew-symmetric (changes sign in the case of odd permutation of variables and does not change under an even permutation of variables) and satisfies the Filippov-Jacobi (sometimes called fundamental identity)
\begin{eqnarray}
[x_1,x_2,\!\!\!&\ldots&\!\!\!,x_{n-1},[y_1,y_2,\ldots,y_n]] =\nonumber\\
    &&\sum_{k=1}^n[y_1,y_2,\ldots,y_{k-1},[x_1,x_2,\ldots,x_{n-1},y_k],y_{k+1},\ldots,y_n].
   \label{n-ary Filippov identity}
\end{eqnarray}
A notion of \emph{$n$-Lie algebra} was independently introduced and studied by V.T. Filippov \cite{Filippov} and Y. Nambu \cite{Nambu}.
In particular case $n=3$ it follows from the definition of $n$-Lie algebra that a \emph{3-Lie algebra} is a vector space $\frak g$ endowed with a \emph{ternary Lie bracket} $[\;,\;,\;]:\frak g\times\frak g\times\frak g\to \frak g$, which has the symmetries
\begin{eqnarray}
[a,b,c]=[b,c,a]=[c,a,b],\;\;\;[a,b,c]=-[b,a,c]=-[a,c,b]=-[c,b,a],
\label{symmetries of ternary commutator 2}
\end{eqnarray}
and satisfies the \emph{ternary Filippov-Jacobi identity}
\begin{eqnarray}
 [a,b[c,d,e]] = [[a,b,c],d,e]+[c,[a,b,d],e]+[c,d,[a,b,e]].
\label{ternary Filippov identity}
\end{eqnarray}
A \emph{metric 3-Lie algebra} is a 3-Lie algebra equipped with a positive definite non-degenerate bilinear form $<\;,\;>$ which satisfies
\begin{equation}
<[a,b,c],d>+<c,[a.b.d]>=0.
\end{equation}
The well known example of simple $n$-Lie algebra can be constructed by means of analog of cross product of vectors of $(n+1)$-dimensional vector space. It was proved in \cite{Ling} that this $n$-Lie algebra is the only simple finite-dimensional $n$-Lie algebra with non-degenerate bilinear form.

\vskip.3cm
\noindent
Given the matrix Lie algebra $\frak{gl}(n)$ of $n$th order matrices one can introduce the ternary Lie bracket of $n$th order matrices as follows
\begin{equation}
[A,B,C]=\mbox{Tr}\,A\;[B,C]+\mbox{Tr}\,B\;[C,A]+\mbox{Tr}\,C\;[A,B].
\label{triple commutator of square matrices}
\end{equation}
Evidently this ternary Lie bracket obeys the symmetries (\ref{symmetries of ternary commutator 2}). It is proved in \cite{Awata-Li-Minic-Yaneya} that this ternary Lie bracket satisfies the ternary Filippov-Jacobi identity and the authors propose to call this ternary Lie bracket \emph{quantum Nambu bracket}. Hence the vector space of $n$th order matrices equipped with the quantum Nambu bracket (\ref{triple commutator of square matrices}) becomes the matrix 3-Lie algebra. The matrix 3-Lie algebra of $n$th order matrices with the triple commutator (\ref{triple commutator of square matrices}) is the metric 3-Lie algebra if we endow it with the bilinear form
$$
<A,B>=\mbox{Tr}\,(AB).
$$

\vskip.3cm
\noindent
Now our aim is to extend this approach to a case more general than the trace of a matrix. Let $\frak g$ be a finite dimensional Lie algebra and $\frak g^\ast$ be its dual space. Fix an element of the dual space $\omega\in \frak g^\ast$ and by analogy with (\ref{triple commutator of square matrices}) define the triple product as follows
\begin{equation}
[x,y,z]=\omega(x)\,[y,z]+\omega(y)\,[z,x]+\omega(z)\,[x,y], \quad x,y,z\in \frak g.
\label{ternary Lie bracket with omega}
\end{equation}
Obviously this triple product is symmetric with respect to cyclic permutations of $x,y,z$ and skew-symmetric with respect to non-cyclic permutations, i.e. it obeys the symmetries (\ref{symmetries of ternary commutator 2}). Straightforward computation of the left hand side and the right hand side of the Filippov-Jacobi identity
\begin{eqnarray}
 [x,y[u,v,t]] = [[x,y,u],v,t]+[u,[x,y,v],t]+[u,v,[x,y,t]],
\label{ternary Filippov identity 1}
\end{eqnarray}
shows that one can split all the terms into three groups, where two of them vanish because of binary Jacobi identity and skew-symmetry of binary Lie bracket. The third group can be split into subgroups of terms, where each subgroup has the same structure and it is determined by one of the commutators $[x,y],[u,v],[u,t],[v,t]$. For instance, if we collect all the terms containing the commutator $[x,y]$ then we get the expression
\begin{equation}
(\omega(u)\,\omega([v,t])+\omega(v)\,\omega([t,u])+\omega(t)\,\omega([u,v]))\,[x,y].
\label{sum of three terms}
\end{equation}
Hence the triple product (\ref{ternary Lie bracket with omega}) will satisfy the ternary Filippov-Jacobi identity if for any elements $u,v,t\in \frak g$ we require
\begin{equation}
\omega(u)\,\omega([v,t])+\omega(v)\,\omega([t,u])+\omega(t)\,\omega([u,v])=0.
\end{equation}
Now we consider $\omega$ as a $\mathbb C$-valued cochain of degree one of the Chevalley-Eilenberg complex of a Lie algebra $\frak g$. Making use of the coboundary operator $\delta:\wedge^k\frak g^\ast\to\wedge^{k+1}\frak g^\ast$ we obtain the degree two cochain $\delta\omega$, where $\delta\omega(x,y)=\omega([x,y])$. Finally we can form the wedge product of two cochains $\omega\wedge\delta\omega$, which is the cochain of degree three and
$$
\omega\wedge\delta\omega(u,v,t)=\omega(u)\,\omega([v,t])+\omega(v)\,\omega([t,u])+\omega(t)\,\omega([u,v]).
$$
Hence the third group of terms (\ref{sum of three terms}) of the Filippov-Jacobi identity vanishes if $\omega\in \frak g^\ast$ satisfies
\begin{equation}
\omega\wedge\delta\omega=0.
\label{condition for omega}
\end{equation}
Thus if an 1-cochain $\omega$ satisfies the equation (\ref{condition for omega}) then the triple product (\ref{ternary Lie bracket with omega}) is the ternary Lie bracket and by analogy with (\ref{triple commutator of square matrices}) we will call this ternary Lie bracket \emph{the quantum Nambu bracket induced by a 1-cochain}.

\vskip.3cm
\noindent
We can generalize the quantum Nambu bracket (\ref{ternary Lie bracket with omega}) if we consider a $\mathbb C$-valued cochain of degree $n-2$, i.e. $\omega\in\wedge^{(n-2)}\frak g^\ast$, and define the $n$-ary product as follows
\begin{equation}
[x_1,x_2,\ldots,x_n]=\sum_{i<j}\,(-1)^{i+j+1}\;\omega(x_1,x_2,\ldots,\hat x_i,\ldots,\hat x_j,\ldots,x_n)\;[x_i,x_j].
\label{n-ary quantum Nambu bracket}
\end{equation}
\begin{thm}
Let $\frak g$ be a finite dimensional Lie algebra, $\frak g^\ast$ be its dual and $\omega$ be a cochain of degree $n-2$, i.e. $\omega\in\wedge^{n-2}\frak g^\ast$. The vector space of Lie algebra $\frak g$ equipped with the $n$-ary product (\ref{n-ary quantum Nambu bracket})
is the $n$-Lie algebra, i.e. the $n$-ary product (\ref{n-ary quantum Nambu bracket}) satisfies the Filippov-Jacobi identity (\ref{n-ary Filippov identity}) if and only if a $(n-2)$-cochain $\omega$ satisfies
\begin{equation}
\omega\wedge\delta\omega=0.
\label{condition for omega in the theorem}
\end{equation}
Particularly the vector space of a Lie algebra $\frak g$ endowed with the $n$-ary Lie bracket (\ref{n-ary quantum Nambu bracket}) is the $n$-Lie algebra if $\omega$ is an $(n-2)$-cocycle, i.e. $\delta\omega=0$.
\label{theorem 1}
\end{thm}
\noindent
Particularly from this theorem it immediately follows that the quantum Nambu bracket for matrices  (\ref{triple commutator of square matrices}) satisfies the ternary Filippov-Jacobi identity. Indeed in this case $\omega=\mbox{Tr}$ and this is 1-cocycle because for any two matrices $A,B$ we have
$$\delta\mbox{Tr}\,(A,B)=\mbox{Tr}\,([A,B])=0.$$
In \cite{Arnlind-Makhlouf-Silvestrov} it is proposed to call the 3-Lie algebra constructed by means of an analog of trace the 3-Lie algebra induced by a Lie algebra. In our approach we will use the similar terminology and call the $n$-Lie algebra constructed by means of a $(n-2)$-cochain $\omega$ (\ref{n-ary quantum Nambu bracket}) the \emph{$n$-Lie algebra induced by a Lie algebra $\frak g$ and an $(n-2)$-cochain} $\omega.$ The $n$-ary Lie bracket will be referred to as the \emph{$n$-ary quantum Nambu bracket induced by a $(n-2)$-cochain}.

\section{Weil algebra of induced 3-Lie algebra}
The aim of this section is to show how one can extend the Weil algebra of a Lie algebra $\frak g$ to the Weil algebra of 3-Lie algebra induced by a Lie algebra $\frak g$.

\vskip.3cm
\noindent
Let $\frak g$ be a finite dimensional Lie algebra and $\frak g^\ast$ be its dual. Let us denote by $\wedge(\frak g^\ast)$ the exterior algebra of $\frak g^\ast$ and by $S(\frak g^\ast)$ the symmetric algebra of $\frak g^\ast$. We remind that the \emph{Weil algebra} of a Lie algebra $\frak g$ \cite{Quillen-Mathai} is the tensor product $S(\frak g^\ast)\otimes\,\wedge(\frak g^\ast)$, which will be denoted by $W(\frak g)$.

\vskip.3cm
\noindent
Let $\mbox{dim}\;\frak g=n$ and $T_a$ be a basis for $\frak g$, where $a$ runs from 1 to $n$. Then
$$
[T_b,T_c]=f^a_{bc}\;T_a,
$$
where $f^a_{bc}$ are the structure constants of a Lie algebra $\frak g$. The structure constants of a Lie algebra satisfies the Jacobi identity
\begin{equation}
f_{bc}^d\;f^k_{ad}+f_{ca}^d\;f^k_{bd}+f_{ab}^d\;f^k_{cd}=0.
\end{equation}
Let $T^{a}$ be the dual basis for $\frak g^\ast$, i.e. $T^a (T_b)=\delta^a_b$. We denote the image of $T^a$ in the exterior algebra $\wedge(\frak g^\ast)$ by $A^a$ and the image of $T^a$ in the symmetric algebra $S(\frak g^\ast)$ by $F^a$. Hence $A^a$ are the generators of the exterior algebra $\wedge(\frak g^\ast)$ and we will assign degree 1 to each of these generators. Analogously $F^a$ are the generators of symmetric algebra $S(\frak g^\ast)$ and we will assign degree 2 to each of these generators. The degree of a generator $A^a$ will be denoted by $|A^a|$ and the degree of $F^a$ will be denoted by $|F^a|$. Then
\begin{eqnarray}
A^a\,A^b \!\!\!&=&\!\!\!(-1)^{|A^a||A^b|}\;A^b\,A^a,\nonumber\\
     A^a\,F^b\!\!\!&=&\!\!\!(-1)^{|A^a||F^b|}\;F^b\,A^a,\nonumber\\
         F^a\,F^b\!\!\!&=&\!\!\!(-1)^{|F^a||F^b|}\;F^b\,F^a.\nonumber
\end{eqnarray}
Now the Weil algebra $W(\frak g)$ is the graded algebra that is $W(\frak g)=\oplus_{i}W^i(\frak g)$. The differential $d:W^i(\frak g)\to W^{i+1}(\frak g)$ is defined by
\begin{equation}
dA^a=F^a-\frac{1}{2}\,f^a_{bc}\,A^b\,A^c,\;\;
       dF^a=f^a_{bc}\,F^b\,A^c.
\label{differential of Weil algebra}
\end{equation}
The differential $d$ has the properties:
\begin{enumerate}
\item
$d$ satisfies the graded Leibniz rule which means that when we differentiate the product of two elements of the Weil algebra and interchange the differential $d$ with the first factor of this product then there appears the coefficient $(-1)^k$ in the second term of Leibniz formula, where $k$ is the degree of the first factor,
\item
$d^2=0$.
\end{enumerate}
It is worth to mention that the Weil algebra is very closely related to the formalism of equivariant differential forms on a principle fiber bundle. Bearing this relation to equivariant forms in mind we can interpret the degree 1 generators $A^a$ of exterior algebra $\wedge(\frak g^\ast)$ as components of a Lie algebra valued connection 1-form on a principal fiber bundle, the degree 2 generators of the symmetric algebra $S(\frak g^\ast)$ as components of a Lie algebra valued curvature 2-form of connection $A^a$ and $d$ as the exterior differential for differential forms.

\vskip.3cm
\noindent
Now our aim is to extend the Weil algebra of a Lie algebra $\frak g$ to the 3-Lie algebra $\triple$ induced by a Lie algebra $\frak g$ with the help of an analog of trace. Let $\omega\in \frak g^\ast$ be an element of the dual space which for any two elements $x,y\in \frak g$ satisfies $\omega([x,y])=0$. Making use of the dual basis for $\frak g^\ast$ we can write $\omega=\omega_a\,T^a$, where the coefficients $\omega_a$ are numbers. Then the condition $\omega([x,y])=0$ is equivalent to the relations
\begin{equation}
f^a_{bc}\;\omega_a=0.
\label{f omega=0}
\end{equation}
The triple commutator of the induced 3-Lie algebra $\triple$ is defined by
\begin{equation}
[x,y,z]=\omega(x)\;[y,z]+\omega(y)\;[z,x]+\omega(z)\;[x,y].
\label{triple commutator section Weil}
\end{equation}
It is shown in the previous section that the triple commutator (\ref{triple commutator section Weil}) satisfies the Filippov-Jacobi identity. Let us denote the structure constants of the induced 3-Lie algebra $\triple$ in a basis $T_a$ by $K^a_{bcd}$ that is
$$
[T_b,T_c,T_d]=K^a_{bcd}\;T_a.
$$
Since the 3-Lie algebra $\triple$ is induced by a Lie algebra $\frak g$ the structure constants of the 3-Lie algebra $\triple$ can be expressed in terms of the structure constants of a Lie algebra $\frak g$ as follows
\begin{equation}
K^d_{abc}=\omega_a\,f^d_{bc}+\omega_b\,f^d_{ca}+\omega_c\,f^d_{ab}.
\end{equation}
It can be proved that the Filippov-Jacobi identity
\begin{equation}
K^a_{bcd}\,K^l_{kma}=K^a_{kmb}\,K^l_{acd}+K^a_{kmc}\,K^l_{bad}+K^a_{kmd}\,K^l_{bca},
\end{equation}
holds in this particular case of induced 3-Lie algebra \cite{AbramovIJGMMP}.
Now we introduce two elements of the Weil algebra $W(\frak g)$ as follows
$$
\chi=\omega_a\,A^a,\; \varphi=\omega_a\,F^a.
$$
Evidently the element $\chi$ has the degree 1 and the element $\varphi$ has the degree 2. If we apply the differential (\ref{differential of Weil algebra}) to the element $\chi$ then we get
\begin{eqnarray}
d\chi\!\!&=&\!\!\omega_a\,d(A^a)=\omega_a\,(F^a-\frac{1}{2}\,f^a_{bc}\,A^bA^c)\nonumber\\
      \!\!&=&\!\!\omega_a\,F^a-\frac{1}{2}\,\omega_a\,f^a_{bc}\,A^bA^c=\omega_a\,F^a=\varphi,
\label{dA omega=F omega}
\end{eqnarray}
where we used the relation (\ref{f omega=0}). Similarly we find that the differential of $\varphi$ is
\begin{equation}
d\varphi=\omega_a\,dF^a=\omega_a\,f^a_{bc}\,F^b\,A^c=0,
\end{equation}
where we again used the relation (\ref{f omega=0}). Hence the element $\varphi$ is closed with respect to the differential $d$ and we can write
\begin{equation}
d\chi=\varphi,\;\;d\varphi=0.
\end{equation}
Let us introduce the following products of elements $A^a,F^a$ and $\chi, \varphi$
\begin{eqnarray}
\chi^a\!\!&=&\!\! \chi\,A^a,\;\;\varphi^a=\varphi\,A^a,\\
\psi^a\!\!&=&\!\! \chi\,F^a,\;\;\Omega^a=\varphi\,F^a.
\end{eqnarray}
These products have the following degrees
$$
|\chi^a|=2,\;\;|\varphi^a|=|\psi^a|=3,\;\;|\Omega^a|=4.
$$
Straightforward computation of the differential of these products gives
\begin{eqnarray}
d\chi^a \!\!\!&=&\!\!\! \varphi^a - \psi^a + \frac{1}{3!}\;K^a_{bcd}\,A^b\,A^c\,A^d,\label{differential A_omega}\\
d\varphi^a\!\!\!&=&\!\!\! \Omega^a - \frac{1}{2}\,f^a_{bc}\,\varphi^a\,A^c,\\
d\psi^a\!\!\!&=&\!\!\! \Omega^a-f^a_{bc}\,\psi^a\,A^c,\\
d \Omega^a\!\!\!&=&\!\!\! f^a_{bc}\;\Omega^b\,A^c.\\
\end{eqnarray}
In order to show more clearly the relation with the structure of the induced 3-Lie algebra $\triple$, we will do the transformation in the space of elements of degree 3 by introducing new elements $\Xi^a_\omega, \Phi^a_\omega$ as follows
\begin{equation}
\Xi^a=\varphi^a+\psi^a,\;\;\Phi^a=\varphi^a-\psi^a.
\end{equation}
Then
\begin{eqnarray}
d\Phi^a\!\!\!&=&\!\!\!\frac{1}{2}\;f^a_{bc}\,A^b\,\Phi^c - \frac{3}{2}\,K^a_{bcd}\,F^b\,A^c\,A^d,\\
d\Xi^a\!\!\!&=&\!\!\!2\,\Omega^a-\frac{1}{2}\,f^a_{bc}\;A^b\,\Xi^c-\frac{1}{2}\,f^a_{bc}\;A^b\,\Phi^c + \frac{1}{2}\,K^a_{bcd}\,F^b\,A^c\,A^d,\\
d\Omega^a\!\!\!&=&\!\!\!K^a_{bcd}\,F^b\,F^c\,F^d.
\end{eqnarray}
Thus in order to construct the extension of Weil algebra $W(\frak g)=S(\frak g^\ast_F)\otimes \wedge(\frak g^\ast_A)$ of a Lie algebra of $\frak g$ to the Weil algebra $W({\triple})$ of the induced 3-Lie algebra $\triple$ we introduce in addition to elements $A^a,F^a$ the element $\chi$ of degree 2, two elements $\Phi^a, \Xi^a$ of degree 3 and the element $\Omega^a$ of degree 4, i.e. we take
\begin{equation}
W(\triple)=S(\frak g^\ast_{F})\otimes S(\frak g^\ast_{\chi})\otimes S(\frak g^\ast_{\Omega})\otimes \wedge(\frak g^\ast_{A})\otimes \wedge(\frak g^\ast_{\Xi})\otimes \wedge(\frak g^\ast_{\Phi}),
\end{equation}
and define the differential by
\begin{eqnarray}
dA^a \!\!\!&=&\!\!\! F^a-\frac{1}{2}\,f^a_{bc}\,A^b\,A^c,\label{1-formula}\\
dF^a \!\!\!&=&\!\!\! f^a_{bc}\,F^b\,A^c,\\
d\chi^a\!\!\!&=&\!\!\! \Phi^a+ \frac{1}{3!}\;K^a_{bcd}\,A^b\,A^c\,A^d,\label{3-formula}\\
d\Phi^a\!\!\!&=&\!\!\!\frac{1}{2}\;f^a_{bc}\,A^b\,\Phi^c - \frac{3}{2}\,K^a_{bcd}\,F^b\,A^c\,A^d,\\
d\Xi^a\!\!\!&=&\!\!\!2\,\Omega^a-\frac{1}{2}\,f^a_{bc}\;A^b\,\Xi-\frac{1}{2}\,f^a_{bc}\;A^b\,\Phi^c + \frac{1}{2}\,K^a_{bcd}\,F^b\,A^c\,A^d,\\
d\Omega^a\!\!\!&=&\!\!\!K^a_{bcd}\,F^b\,F^c\,F^d.\label{6-formula}
\end{eqnarray}
It can be proved that the differential $d$ defined by the above formulae satisfies $d^2=0$. Hence the formulae (\ref{3-formula}) -- (\ref{6-formula}) can be regarded as the extension of differential $d$ defined on the Weil algebra $W(\frak g)$ to differential on the Weil algebra $W(\triple)$ of the induced 3-Lie algebra.

\vskip.3cm
\noindent
The formulae can be written by means of commutators (binary and ternary) if we consider the tensor product $\frak g\otimes W(\triple)$. Given three elements of this product $X_1\otimes\Lambda_1, X_2\otimes\Lambda_2, X_3\otimes\Lambda_3$, where $X_1,X_2,X_3\in\frak g$ and $\Lambda_1,\Lambda_2,\Lambda_3\in W(\triple)$ we can compose two products as follows
\begin{eqnarray}
&&[X_1\otimes\Lambda_1, X_2\otimes\Lambda_2]=[X_1,X_2]\otimes \Lambda_1\cdot\Lambda_2\\
&&[X_1\otimes\Lambda_1, X_2\otimes\Lambda_2, X_3\otimes\Lambda_3]=[X_1,X_2,X_3]\otimes\Lambda_1\cdot\Lambda_2\cdot\Lambda_3.
\end{eqnarray}
Let us denote 
\begin{eqnarray}
A &=& T_a\otimes A^a,\; F=T_a\otimes F^a,\; \chi=T_a\otimes \chi^a,\nonumber\\
\Phi &=& T_a\otimes \Phi^a,\; \Xi=T_a\otimes \Xi^a,\; \Theta=T_a\otimes\Theta^a.\nonumber
\end{eqnarray}
Then the formulae (\ref{1-formula})--(\ref{6-formula}) can be written in the form
\begin{eqnarray}
dA\!\!\!&=&\!\!\! F-\frac{1}{2}\,[A,A], \nonumber\\
dF\!\!\!&=&\!\!\! [F,A],\nonumber\\
d\chi\!\!\!&=&\!\!\! \Phi+ \frac{1}{3!}\;[A,A,A],\nonumber\\
d\Phi\!\!\!&=&\!\!\!\frac{1}{2}\,[A,\Phi] - \frac{3}{2}\,[F,A,A],\nonumber\\
d\Xi\!\!\!&=&\!\!\!2\,\Omega-\frac{1}{2}\,[A,\Xi]-\frac{1}{2}[A,\Phi] + \frac{1}{2}\,[F,A,A],\nonumber\\
d\Omega\!\!\!&=&\!\!\![F,F,F].\nonumber
\end{eqnarray}
\section{B.R.S. algebra and 3-Lie algebra}
Our aim in this section is to show how one can use the 3-Lie algebra $\triple$ to construct the universal B.R.S. algebra \cite{Dubois-Violette_Talon_Vialett}. The advantage of approach based on the 3-Lie algebra $\triple$ is that ghost field, auxiliary field and their transformations naturally arise from the structure of 3-Lie algebra.

\vskip.3cm
\noindent
First we remind the notion of universal B.R.S. algebra associated to a Lie algebra $\frak g$. Let $T_a$ be a basis for $\frak g$ and $T^a$ be the dual basis for $\frak g^\ast$. Take four copies of the dual space $\frak g^\ast_A,\frak g^\ast_F,\frak g^\ast_\chi,\frak g^\ast_\phi$ with dual basis $A^a,F^a,\chi^a,\phi^a$ respectively. Let $\mathscr A(\frak g)$ be the free bigraded commutative algebra generated by $A^a$ in the bidegree $(1,0)$, $\chi^a$ in the bidegree $(0,1)$, $F^a$ in the bidegree $(2,0)$, $\phi^a$ in the bidegree $(1,1)$, i.e.
$$
\mathscr A(\frak g)=\wedge(\frak g^\ast_A)\otimes S(\frak g^\ast_F)\otimes\wedge(\frak g^\ast_\chi)\otimes S(\frak g^\ast_\phi).
$$
In analogy with the previous section we introduce the tensor product $\frak g\otimes\mathscr A(\frak g)$ and the products of elements $A = T_a\otimes A^a,\; F=T_a\otimes F^a,\; \chi=T_a\otimes \chi^a,\phi=T_a\otimes\phi^a$. Define two differentials $d,\delta$, of bidegrees $(1,0)$ and $(0,1)$ respectively as follows
\begin{equation}
dA=F-\frac{1}{2}\;[A,A], \;dF=[F,A], \;d\chi=\phi,\;d\phi=0,
\end{equation}.
and 
\begin{equation}
\delta A=-\phi-[A,\chi],\;\delta F=[F,\chi],\;\delta\chi=-\frac{1}{2}\;[\chi,\chi],\;\delta\phi=[\phi,\chi].
\end{equation}
These differentials sarisfy the equations
$$
d^2=\delta^2=(d+\delta)^2=0.
$$
The algebra $\mathscr A(\frak g)$ is referred to as the universal B.R.S. algebra of a Lie algebra $\frak g$.

\vskip.3cm
\noindent
Now our aim is to use the structure of the 3-Lie algebra $\triple$ to construct the universal B.R.S. algebra. For this purpose we use the idea of gauge transformation by means of triple commutator
\begin{equation}
[X,Y,Z]=\mbox{Tr}\,X\;[Y,Z]+\mbox{Tr}\,Y\;[Z,X]+\mbox{Tr}\,Z\;[X,Y],
\label{formula 1 in introduction A}
\end{equation}
proposed in\cite{Awata-Li-Minic-Yaneya}. Indeed one can consider the following analog of gauge transformation
\begin{equation}
D\,A=i\,[X,Y,A].
\end{equation}
Following this idea we introduce two elements $\xi,\eta$ of the bidegrees $(0,1),(0,0)$ respectively and define the differential $\delta$ of $A$ as follows
\begin{equation}
\delta A=-[A,\xi,\eta].
\label{gauge transformation}
\end{equation}
Applying the formula for the triple Nambu bracket we find
\begin{equation}
\delta A=-[\mbox{Tr}(\xi)\,\eta-\mbox{Tr}(\eta)\,\xi, A]-\mbox{Tr}(A)\;[\xi,\eta].
\label{delta A}
\end{equation}
Now it is natural to define the ghost field $\chi$ and the auxiliary field $\phi$ by
\begin{equation}
\chi=\mbox{Tr}(\xi)\,\eta-\mbox{Tr}(\eta)\,\xi, A, \;\phi=\mbox{Tr}(A)\;[\xi,\eta]
\label{chi phi}
\end{equation}
Substituting these definitions into (\ref{delta A}) we obtain the differential for $A$, i.e.
$$
\delta A=-\phi-[A,\chi].
$$  
Now we define the differential $\delta$ for the elements $\chi,\phi$ by
\begin{eqnarray}
\delta\xi \!\!&=&\!\!-\frac{1}{2}\;\mbox{Tr}(\eta)\,[\xi,\xi]+\mbox{Tr}(\xi)\,[\eta,\xi],\label{delta xi}\\
\delta\eta\!\!&=&\!\!-\mbox{Tr}(\eta)\,[\eta,\xi].\label{delta eta}
\end{eqnarray}
Differentiating by $\delta$ the both sides of (\ref{chi phi}) and making use of (\ref{delta xi}), (\ref{delta eta}) we get the formulae for the differential $\delta$ of $\chi$ and $\phi$, i.e.
$$
\delta\chi=-\frac{1}{2}\;[\chi,\chi],\;\delta\phi=[\phi,\chi].
$$
We conclude this section by pointing out that proposed here approach clearly demonstrates that the formulae for the differential $\delta$ (or BRS-operator) are consequences of two-parameter gauge transformation (\ref{gauge transformation}) constructed with the help of triple quantum Nambu commutator (\ref{formula 1 in introduction A})
\subsection*{Acknowledgment}
The author gratefully acknowledges that this work was financially supported by the institutional funding IUT20-57 of the Estonian Ministry of Education and Research.


\begin{thebibliography}{1}
\bibitem{AbramovIJGMMP}
Abramov, V., \textit{Matrix 3-Lie superalgebras and BRST supersymmetry}, to appear in Int. J. of Geom. Meth. in Mod. Phys, {\bf 14} (2017),  arXiv:1707.01270.
\bibitem{Abramov2017Gent}
Abramov, V., \textit{Quantum super Nambu bracket of cubic supermatrices and 3-Lie superalgebra}, to appear in the Proceedings of ICCA conference in Gent.
\bibitem{Awata-Li-Minic-Yaneya}
H. Awata, M. Li, D. Minic, and T. Yaneya, \textit{On the quantization of Nambu brackets}, {JHEP02} (2001) 013.
\bibitem{Dubois-Violette_Talon_Vialett}
Dubois-Violette, M., Talon, M., Vialette, C.M., \textit{B.R.S. Algebras. Analysis of Consistency Equations in Gauge Theory}, Commun. Math. Phys. {\bf 102} (1985), 105--122.
\bibitem{Filippov}
V.T. Filippov, \textit{$n$-Lie algebras}, {Siberian Math. J.} {\bf 26} (1985), 879--891.
\bibitem{Ling}
Ling, W.X., \textit{On the structure of $n$-Lie algebras}, PhD thesis, Siegen, 1993.
\bibitem{Nambu}
Y. Nambu, \textit{Generalized Hamiltonian mechanics}, {Phys. Rev.} D, {\bf 7} (1973), 2405--2412.
\bibitem{Quillen-Mathai}
Mathai, V., Quillen, D., {"Superconnections, Thom Classes, and Equivariant Differential Forms"}, Topology \textbf{25}, 1, 85--110 (1986).
\bibitem{Takhtajan}
L. Takhtajan, {On foundation of generalized Nambu mechanics}, {\it Comm. Math. Phys.} {\bf 160} (2) (1994), 295 -- 315.
\end{thebibliography}
\end{document}